\documentclass{amsart}

\usepackage{amsmath, amsthm}
\usepackage{amssymb, latexsym}
\usepackage{amsfonts}

\newtheorem{theorem}[equation]{Theorem}

\newtheorem{proposition}[equation]{Proposition}

\newtheorem{Example}[equation]{Example}

\newtheorem{Remark}[equation]{Remark}

\theoremstyle{definition}

\newcommand{\beql}[1]{\begin{equation}\label{#1}}
\newcommand{\eeq} {\end{equation}}

\font\Bbb=msbm10

\newcommand\R{\hbox{\Bbb R}}
\newcommand\Z{\hbox{\Bbb Z}}

\font\Bbb=msbm10

\def\R{\hbox{\Bbb R}}

\def\Z{\hbox{\Bbb Z}}

 \def\CC{\hbox{\Bbb C}}

 \font\Bbcc=msbm7
 
 \def\CCCC{\hbox{\Bbcc C}}

\numberwithin{equation}{section}

\DeclareMathOperator{\Aut}{Aut}

\newcommand\e{{ \epsilon }}

\newcommand\G{\Gamma }
\newcommand\D{{ \Delta }}
\newcommand\z{{ \zeta }}

        \def\GF{{\rm GF}}

        \def\Sz{{\rm Sz}}

        \def\<{{\langle}}
        \def\>{{\rangle}}

\font\Bbb=msbm10

\def\R{\hbox{\Bbb R}}

\def\Z{\hbox{\Bbb Z}}

\font\Bbbsmall=msbm7

\def\Zsmall{\hbox{\Bbbsmall Z}}

\def\a{\alpha}

\def\div{ \kern-.5pt\hbox{\big |} }
\def\ndiv{ {\not\kern-.5pt\hbox{\big |}\,} }
\def\ndivv{ {\not\kern+1.5pt\hbox{$\mid$}\,} }

\def\Aut{{\rm Aut}}

\def\dim{{\rm dim}}

\def\B{^2\kern-.8pt B}
\def\G{^2\kern-.8pt G}
\def\EH{^2\kern-.8pt\hat  E}
\def\E{^2\kern-.8pt E}
\def\D{^3\kern-1pt D}
\def\FF{^2\kern-.8pt F}

\newdimen\refcodesize
\newbox\seriesbox
\def\GF{{\rm GF}}

\def\ZZ{{\Z}}

\newcommand{\sF}{{\mathcal F}}
\newcommand{\zen}{{\mathcal Z}}

\newcommand{\oA}{\overline{A}}

\newcommand{\oB}{\overline{B}}

\DeclareRobustCommand{\SkipTocEntry}[4]{}
\begin{document}

\title[MUBs inequivalence  and affine planes]
{MUBs inequivalence  and affine planes}
\thanks{This research was supported in part by
 NSF grant DMS~0753640.}
     
       \author{W. M. Kantor}
       \address{Department of Mathematics, University of Oregon,
       Eugene, OR 97403 }
             \email{kantor@uoregon.edu}

\begin{abstract}

There are fairly large families of unitarily inequivalent complete sets of  $N+1$ mutually unbiased bases (MUBs) 
known
in 
$\CCCC^N$ for various prime powers $N$.
The number of such sets is not bounded above by any polynomial as a function of $N$.
 While it is standard that there is a superficial similarity  
 between  complete sets of  MUBs
 and  finite affine planes,
 there is an intimate relationship between  these large families  and affine planes. This note  briefly summarizes ``old'' results 
   that do not appear to be well-known
   concerning known families of complete sets of  MUBs  and their associated planes. 
   
\vspace{-8pt}

\end{abstract}

\maketitle
\centerline{\em In memory of  Jaap Seidel}

\vspace{-2pt}
  \section{Introduction} 
\label{Introduction}

Starting with \cite{Al, Iv, W, WF},
there has been a great deal of activity constructing, studying and using ``complete sets of mutually unbiased bases'' (MUBs)
of $\CC^N$, which are known to exist when $N$ is a prime power (see  \cite{GR} and the references therein).
It   is proved in \cite{GR}  that almost all published constructions 
produce unitarily equivalent complete sets of MUBs, 
while it is also observed that there are other complete sets. 
The purpose of this note is to indicate  the large number of 
{other} known complete sets of MUBs,
together with hints of the geometric context of those constructions and non-equivalences.  Almost all of these known sets fall into the same framework (Theorem~\ref{Char 2 MUBs} for arbitrary prime powers); the only known exceptions are in Example~\ref{planar function examples}.

Whereas sets of MUBs are usually viewed in terms of  sets
of vectors, it is easier to discuss automorphisms and equivalence 
(defined in Section~\ref{MUBS characteristic 2}) by using the 1-spaces spanned by those vectors%
\footnote{We do not use bases since automorphisms   
do not preserve   bases (e.\hspace{.2pt}g., $Z(b)$
 in (\ref{X and Z})
does not preserve the standard basis).}.  
This leads to the following definitions.

An {\em orthoframe}\footnote{We avoid the term {\em frame} used in
 \cite{CCKS} so as not to conflict with other   uses  for that word.}
   in $\CC^N$ is a family of $N$ pairwise orthogonal 1-spaces. Distinct  orthoframes $\sF_1, \sF_2$ are called 
\emph{mutually unbiased} if 
$| (u_1,u_2) | =  1/\sqrt N$ whenever $u_i$ is a unit vector of a member  of $\sF_i$ for $i=1,2$.  If a unit vector is taken from each member of an orthoframe the result is an orthonormal basis;   orthonormal bases are called {\em mutually unbiased} (MUBs) if the same is true for the corresponding orthoframes.  
{\em We will   view mutually unbiased bases and mutually unbiased orthoframes as}  ``essentially'' {\em the same objects.} 
Any family of MUBs has size at most $N+1$.  A set 
$\sF$
of MUBs meeting this bound is called {\em complete} or {\em maximal}.
 
The  construction of the  sets $\sF$ highlighted here  occurred in the early 1990's.    
The authors of \cite{CCKS}  did not know 
the 
Physics context 
(in fact they were concerned with
 the union $\cup\sF$ of the members of $\sF$).
Their  construction method was based on 
connections with affine planes%
\footnote{An \emph{affine plane of order $N$} is a 
combinatorial object consisting of a set of 
$N^2$ \emph{points}, together with   $N^2+N$ point-sets of size $N$ called 
\emph{lines}, such that any two distinct points are on a unique line.  Then the lines fall into $N+1$ ``parallel classes'' of size $N$, each of which partitions the points.},
symplectic spreads (defined in Section~\ref{MUBS characteristic 2}),  sets $\sF$
 and error-correcting codes over $\Z_2$ or $\Z_4$.  Symplectic spreads were also observed 
\cite{KaOSU}
to be related to orthogonal decompositions of the Lie algebras ${\rm sl}_n(\CC)$;  this was recently rediscovered in  \cite{BSTW2}.  Relationships between symplectic spreads, $\Z_p$-codes for odd $p$ and unimodular lattices appear in \cite{ST}.
A number of papers, such as \cite{StHP},   use    complete sets of MUBS for CDMA  (code division multiple access)  in radio communication technologies.
The research in \cite{CCKS} was initiated by Seidel's   interest in    sets $\sF$  related to cubature formulas
\cite{DGS, Se, Ko}.

The subject matter of this brief note has been surveyed already
(e.\hspace{1pt}g., in 
\cite{KaAMS}), but emphasizing only planes and  codes, not 
$\CC ^N$ or $\R^N$.
The motivation behind this note is that many of the results in 
\cite{CCKS}  are not widely known to the MUBs community,
although \cite{CCKS} is referenced in several MUBs papers.
In particular, \cite{GR,RS} are essentially  the only recent references that observe  that there are several inequivalent 
complete sets of MUBs for some dimensions $N$.
As a function of $N$ 
 the number of such sets  
is not bounded above by any polynomial  
 (Example~\ref{symplectic char 2}(a)).

The simplest examples (Examples~\ref{2-dim examples}, \ref{odd 2-dim examples}(a)) arise from 2-dimensional vector spaces, and hence from affine planes over finite  fields. 
It is unknown whether or not there is a general relationship between complete sets of MUBs and affine planes.
  All known finite  affine planes
have prime power order; all known complete sets of MUBs have 
$p^n+1$ members for some prime $p$. This may or may not be a 
coincidence\footnote{See the delectable
 observation at the end of \cite{Be}.}, but we will have nothing to say about planes or complete sets of MUBs that do not arise from a prime power.
However, Remark~\ref{Where are the planes} can be interpreted as a somewhat negative observation concerning the occurrence of an affine plane
in a prime power instance.

Sections~\ref{MUBS characteristic 2} and \ref{MUBS  $p>2$} deal with complete sets of MUBs in
$\CC^{p^n}$ for  $p=2$ and $p>2$, respectively. 
  Section~\ref{Real MUBs}   briefly considers the situation when complex spaces are replaced by real ones.  
We have not discussed the quaternionic version of complete sets of MUBs.  As in the real and complex cases,  these  are plentiful \cite{quaternionic}.

\section{complete sets of MUBS in $\CC^{2^n}$}
\label{MUBS characteristic 2}
  
Equip $V=\ZZ_2^n$  with its usual dot product  $x \cdot y$,
 and $\CC^N, N=2^n$, with its usual hermitian inner product
 $(~,~)$.
Label the standard basis of $\CC^N$ as $e_v$,
$v \in V$.
For $b \in V$   define linear transformations
$X(b)$ and  $Z(b)$
 on 
$\CC^N$    
 by 
\begin{equation}\label {X and Z}
\mbox{$X(b) \colon  e_v \mapsto e_{v+b}$ \quad and
\quad  $Z(b): e_v \mapsto (-1)^{b \cdot v} e_v$.}
\end{equation}

The groups $X(V):  = \{X(b) \mid b \in V \}$ and
 $Z(V) :  = \{ Z(b) \mid b \in V \}$ consist of unitary transformations and  are isomorphic to   the additive group $V$. 
Moreover,  the group $E :=  X(V) Z(V) \{\pm I \}$ they generate is an \emph{extraspecial group} (or {\em Heisenberg group}) of order $2^{1+2n}$ with center
 $\zen(E)=  \{\pm I \}$, which we identify with $\ZZ_2$.
We also need the slightly larger group $P:= 
 E \{\pm I, \pm iI \} $  for the usual $i\in \CC$, with center $\zen(P)=\{\pm I, \pm iI \} $ (see the comments preceding Proposition~\ref{MUB test} below).  
We use the natural map ~${}^{\overline{\ \  }}\colon\! P\to \overline{P} =P/\zen(P)\cong V\oplus V$, and therefore avoid using complex conjugation in calculations.
 The commutator 
\begin{equation}
\label{commutator}\big (X (a) Z (b) \big)^{-1} \big(X (a') Z (b') \big)^{-1}
 \big ( X (a) Z (b)\big )\big (X (a') Z (b')\big ) =  
a \cdot b' - a'\cdot b
\end{equation}
on $P$  determines a nondegenerate alternating bilinear form {\bf (~,~)} on the $\Z_2$-space $P/\zen(P)\cong V\oplus V \cong E/\zen(E)$. 

If $A$ is an abelian subgroup of $P$ such that $\overline A$ is a {\em totally 
isotropic} $n$-space   of $\overline{P}$ (i.\hspace{1pt}e., 
$\dim \overline A= n $ and
{\bf ($\overline{A}$,$\overline{A}$)}$=0$), then the set 
  $\sF (A)=\sF (\overline A)$ of $A$-irreducible subspaces of $\CC^{N}$
is an orthoframe. (Equivalently: $\sF (A)$ is the set of 
all 1-dimensional subspaces invariant under   $A$.)  Moreover, $\sF (A)$ is   invariant under $P$. If $B$ is a second such subgroup of $P$ for which 
  $\oA \cap \oB = 0$,   
then $| (u_1,u_2) | =  2^{-n/2}=1/\sqrt N$ whenever $u_1$ and $u_2$ are unit vectors in members of $\sF (A) $ and $ \sF (B)$, respectively
\cite[Lemma~3.3]{CCKS}:
$\sF (A) $ and $ \sF (B)$ are mutually unbiased.

{\em Each totally isotropic $n$-space} of  $\overline{P}$ arises as
some  $\overline{A}$, and hence  {\em determines a unique orthoframe $\sF (\overline{A}) $.}

 A {\em symplectic spread} of the symplectic space  $\overline{P}$ (or   of $\overline{E} =E/\zen(E)$)
 is a family $\Sigma$ of
 $N+1$ totally isotropic $n$-spaces  of  $\overline{P}$   any two of which have  intersection $0$.   Then every nonzero vector of  $\overline{P}$  is in one and  only one member of $\Sigma$
 (so $\Sigma $ partitions the nonzero vectors). {\em This determines an  affine plane of order} $N$, whose points are the vectors in  $\overline{P}$   and whose lines are the translates of the members of $\Sigma$ by the elements of   $\overline{P}.$    This is the elementary relationship between  affine  planes and the sets of MUBs considered in this note  (see \cite{KaAMS} and (\ref{affine plane})).
 
 The rest of this  survey is concerned with the following result and its consequences and variations (such as its validity for odd characteristic).
 
 \begin{theorem} {\rm\cite [Theorem~5.6 and Proposition~5.11]{CCKS}}~
 \label{Char 2 MUBs}
Each symplectic spread  $\Sigma$ of $\overline{P}$ determines a complete set
 $\sF (\Sigma) =\{ \sF (\overline{A}) \mid  \overline{A}\in \Sigma\} $ of $N+1$  {\rm MUBS}  in $\CC^N $ such that each  member is 
 invariant under $P$.

Let $\Sigma'$  be another symplectic spread of $\overline{P}.$ 
Then   $\Sigma $ and $\Sigma'$ are equivalent under a linear transformation  of $\overline{P}$ preserving the alternating bilinear form on $\overline{P}$ if$,$ and only if$,$ 
$\sF (\Sigma)$ and $\sF (\Sigma')$  are equivalent under a unitary transformation of $\CC^N \!$.
$($This occurs if and only if the corresponding affine planes are isomorphic.$)$
 \end{theorem}
  
  Two complete sets of MUBs are called (unitarily) \emph{equivalent} if there is a unitary transformation sending one set to the other as in the theorem.
  
\begin{Example}
\label{2-dim examples}
\ \rm  The 1-dimensional vector space    $V= \GF(2^n)$
over a finite field $ \GF(2^n)$  of size $2^n$ is also an $n$-dimensional vector space over $\GF(2)=\ZZ_2$.
Let $T\colon\GF(2^n)\to \ZZ_2$ be the trace map
(so that $T(x):=\sum_{i=0}^{n-1}x^{2^i}$), and use the alternating bilinear form 
{\bf($(a,b),(c,d)$)}$:=T(ad-bc)$ on  the $2n$-space   $V\oplus V=\GF(2^n)^2$ over $\GF(2)$.  Then the set $\Sigma$ of 1-dimensional $ \GF(2^n)$-spaces of $V\oplus V$ is%
\footnote{We are identifying isomorphic vector spaces.}
 a symplectic spread of 
$\overline{P} =V\oplus V$ (since the determinant $ad-bc$ vanishes on each of them).  This produces a complete set  $\sF (\Sigma)$ of MUBs.
As proved in   \cite{GR},  $\sF (\Sigma)$  is equivalent to the complete sets obtained in most previous papers. 

We emphasize that there is nothing mysterious about $\Sigma$:
it is just the subsets $x=0$ and $y=mx$ of $ \GF(2^n)^2$ for $m\in \GF(2^n)$. This should be reminiscent of the  lines through the origin in  high school
 (cf.~(\ref{affine plane})).
The orthoframes determined by the members of $\Sigma$ are described using   sums involving complex roots of unity, just as in many references such as \cite{BBRV,GR,WF}\, (cf. (\ref{hats}) for explicit 1-spaces). 
\end{Example}

\begin{Example}
\label{old examples}
Some inequivalent examples. \rm
There are many other known symplectic spreads of  $\overline{P}=V\oplus V$ for suitable $n$. They are complicated to describe (this is discussed at length in \cite{KW1}; cf. (\ref{hats})).  We present an example taken from  \cite{KaKer1}.

As in Example~\ref{2-dim examples},
let   $V=\GF(2^n)$  
and equip   the $\GF(2)$-space
 $ V\oplus V$   with  the previous bilinear form.
 Assume that   $n>3$ is odd.
Then \emph{the subsets 
$x=0$ and $y=m^2x+mT(x)+T(mx)$ of $V\oplus V,$ $m\in V,$ are a symplectic spread of $V\oplus V$ that is not equivalent to the one in} Example~\ref{2-dim examples}.

The complete set  of MUBs produced by this example and 
the one in Example~\ref{2-dim examples} are not unitarily equivalent, in view of Theorem~\ref{Char 2 MUBs}.
In order to write explicit vectors in $\CC^{2^n}$ we would need to lift all of this from $ \Z_2$ to $\Z_4$,  as discussed at length in \cite[Sec.~5]{CCKS};
cf.  (\ref{hats}).  (For odd prime powers a  simpler example, including
explicit  complex  vectors,  appears in Example~\ref{odd examples}(b).)

\end{Example}

We now provide a list  indicating that there are many different families with different properties. 
First   we need to discuss one obvious  aspect of any complete set  $\sF$ of MUBs: its automorphism group $\Aut  (\sF)$.  
This consists of all unitary transformations of $\CC^N$ that send 
$\sF$ to itself.
 We have already seen that $P$ 
lies in $\Aut  (\sF)$,
inducing the identity on  $\sF$
(cf.  Proposition~\ref{MUB test});
the same is true  for
 all unitary  matrices $\a I$ with $\a\in \CC$,
 $|\a|=1$.
Using  Theorem~\ref{Char 2 MUBs}, {\em every automorphism of $\sF(\Sigma)$ normalizes $P$ and sends $\Sigma$ to itself.}
Acting by conjugation,
 $\Aut  (\sF)$ then
  induces on $\overline{P}$   a subgroup 
$\overline{\Aut}  (\sF)$
of the symplectic group  
 (in  fact, $\overline{\Aut}  (\sF)$   is the set-stabilizer of $\Sigma$
in the symplectic group of isometries of  $\overline{P}$).  This is
part of 
 what we will focus on when describing some pairwise inequivalent  complete sets of  MUBs.

\begin{Remark}
\label {2-dim examples'}
\rm In Example~\ref{2-dim examples},
 $\overline{\Aut}  (\sF)$   contains  all invertible semilinear transformations 
 $(x,y)\mapsto M(x^\sigma,y^\sigma) $ on $\GF(2^n)^2$
with $\sigma \in \Aut(\GF(2^n))$ and $M $ a 
 linear transformation on the $2$-space $\GF(2^n)^2$ of determinant $1$.   (For all other {\em known}  examples     
$\overline{\Aut}  (\sF)$  is  much smaller.)

For these examples 
${\Aut}  (\sF)$ is 3-transitive on $\sF(\Sigma)$, and has a cyclic subgroup of order $2^n+1$ 
that is transitive on $\sF(\Sigma)$ (cf. 
Example~\ref{symplectic char 2}(c)     below).
There is also a cyclic subgroup of order $2^n-1$ fixing two members of 
   $\sF(\Sigma)$ and permuting the remaining ones transitively
    (cf. 
Example~\ref{symplectic char 2}(b)     below).
   
\end{Remark}
  
\begin{Example} \rm 
\label{symplectic char 2}
We indicate some of the {\em other} known examples of 
complete sets of  MUBs in $\CC^N = \CC^{2^n}$  arising from
symplectic spreads of   $\Z_2^{2n}$, together with  additional remarks concerning them.  
We emphasize that examples (a)-(d) occur in $\CC^{2^n}$ with $n>3$ not a power of $2$;
only in (a) can $n$  be  prime.
(``Unbounded'' means   as a function of $N=2^n$.)
\begin{itemize}
\item[(a)] Examples $\sF(\Sigma)$ in $\CC^{2^n}  $ for which  
\emph{$\overline{\Aut}  (\sF)$  is an extension of the additive group    
$\GF(2^n)^+$ by a subgroup of $\GF(2^n)^*\Aut(\GF(2^n))$},
where ${n>3}$ is   not a power of $2$    \cite{KW2}: ~
$ {\Aut}  (\sF)$  has an elementary abelian subgroup of order $2^n$  that induces the identity on  one member of $ \sF(\Sigma)$ and is transitive on the remaining ones.  

\emph{The number of pairwise unitarily inequivalent   complete sets of} MUBs \emph{of this sort is not bounded above by any polynomial in $N$.}  The number of these complete sets is an increasing function of the number of prime divisors of~$n$.

The smallest $N$ for which there are inequivalent complete sets of MUBs  obtained via Theorem~\ref{Char 2 MUBs} is $N=2^5$.

\smallskip
Below, ``unbounded"  will mean  as a function of $N$.

\item[(b)]  An unbounded number of  examples  \emph{$\sF(\Sigma)$ 
in $\CC^{2^n} \! $  for which  
$\overline{\Aut}  (\sF)$  is an extension of the multiplicative group of  
$\GF(2^n)$ by a subgroup of $\Aut(\GF(2^n)),$}
where $n$ has at least two odd prime factors  \cite{KW3}:
$\Aut  (\sF(\Sigma))$ 
has a cyclic subgroup of order $2^n-1$ fixing
 a pair of members of $ \sF(\Sigma)$ and  transitive on the remaining ones.  (For many of these sets 
 $\sF(\Sigma)$ of MUBs
each member  of the indicated pair of orthoframes is sent to itself
by $\Aut  (\sF(\Sigma))$;
 for others these two are interchanged.) 

\item[(c)]  An unbounded number of  examples  \emph{$\sF(\Sigma)$ 
in $\CC^{2^n} \! $  for which  
$\overline{\Aut}  (\sF)$  is an extension of 
a cyclic group of order $2^n+1$ by a subgroup of $\Aut(\GF(2^n)),$} 
 where $n$ is neither prime nor a power of  $2 $  \cite{KW1}: 
$\Aut  (\sF(\Sigma))$  has a  \emph{cyclic  subgroup  of order $2^n+1$ that is  transitive on the family $ \sF(\Sigma)$},
as in Example~\ref {2-dim examples'}.

\item[(d)] An unbounded number of  examples $\sF(\Sigma)$
in $\CC^{2^n} \! $  \emph{with 
  $n>9$ odd and composite,   and 
$\overline{\Aut}  (\sF) =1$ \cite{KaBoring}}:   there is a great deal of structure available
for these examples, enough to prove that the automorphism group  $\Aut  (\sF(\Sigma))$ is remarkably small and nevertheless to be able to prove   inequivalences.

\item[(e)] There is a symplectic spread $\Sigma$ in $\Z_2^{n}$, 
 $n\equiv 4$ (mod 8),  $n>4$,  arising from the Suzuki group $\Sz(2^{n/4})$  \cite[Prop.~3.3]{Ti}, so that $\sF(\Sigma)$ is a complete set of MUBs in 
 $\CC^{2^n} \! $.

\end{itemize}

Each of these families produces sets of 1-spaces of  $\CC^N$ that can be described explicitly using $\Z_4^n$ together with  $V=\Z_2^n$ 
 \cite[Sec.~5]{CCKS}; cf. (\ref{hats}).

Each of the  known families (b)-(e) contains fewer than   $\sqrt N$ pairwise inequivalent  complete sets of MUBs in $\CC^N$,~which is quite different from the situation in~(a). 

Examples (a)-(d) were obtained using what amounts to an algorithm that starts with Example~\ref{2-dim examples} and uses quadratic and alternating bilinear forms on $\Z_2$-spaces together with field changes  
(\cite[Sec.~3]{KaAMS}, \cite[Secs.~2.6,~2.7]{KW2}).  This approach only works in characteristic 2.  There are undoubtedly large numbers of other examples yet to be found.

\end{Example}

We have focused on $P$.  We could just as well have used the slightly smaller extraspecial group $E$  for the purpose of describing constructions (though not for full automorphism groups or proving inequivalence!). 
Namely, a preimage of a totally isotropic $n$-space of 
$\overline{E}$ is diagonalizable in $\CC^{2^n}\!$ using 
a unique orthoframe.  
Moreover, it is not difficult to use extraspecial groups to test whether or not a given complete set of MUBs arises as in Theorem~\ref{Char 2 MUBs}:
\begin {proposition}
\label{MUB test}
A complete set $\sF$ of {\rm MUBS} in $\CC^{2^n}\!$ arises as in 
{\rm Theorem~\ref{Char 2 MUBs}} if and only if there is an extraspecial group of  $2^{1+2n}$ unitary transformations   sending each member of $\sF$ to itself.
\end {proposition}

A starting point for \cite{CCKS} was the study of ``Kerdock codes'' over $\Z_2$ and $\Z_4$.  In \cite{Ca} the term
``Kerdock codes'' was redefined to 
be sets of  $N^2+N$ unit vectors in  orthonormal bases of $\CC^N$ determined by 
$\cup \sF$ for  $\sF$  in Theorem~\ref{Char 2 MUBs}.

\section{complete sets of MUBS in $\CC^{p^n}\!,$ $p>2$}
\label{MUBS  $p>2$}  

Consider an odd prime $p$ and $V=\ZZ_p^n$  with its usual dot product  $x \cdot y$.
Equip $\CC^N,$ $ N=p^n$, with its usual hermitian inner product
 $(~,~)$.
Label the standard basis of $\CC^N$ as $e_v$,
$v \in V$.
Let $\z\in \CC $ be a primitive $p^{\rm th}$ root of unity.
For $b \in V$,   define  
$$
\mbox{$X(b) \colon  e_v \mapsto e_{v+b}$ \qquad and
\qquad  $Z(b): e_v \mapsto \z^{b \cdot v} e_v $.}
$$
The groups $X(V):  = \{X(b) \mid b \in V \}$ and
 $Z(V) :  = \{ Z(b) \mid b \in V \}$ consist of unitary transformations and  are isomorphic to   the additive group $V$. 
Moreover, they generate   an extraspecial group 
(or Heisenberg group)
 $E :=  X(V) Z(V) \{\z^j I\mid 0\le j<p\} $
of order $p^{1+2n}$ with center
 $\zen(E)= \{\z^j I\mid 0\le j<p\}$, which we identify with $\ZZ_p$.
We use the natural map \ ${}^{\overline{\ \  }}$ \ as  before. 
~
 The commutator (\ref{commutator})  again  defines a nondegenerate alternating bilinear form on $E/\zen(E)\cong V\oplus V$.

If $A$ is an abelian subgroup of $P$ such that $\overline A$ is a totally 
isotropic $n$-space   of $\overline{E}$, then the set 
 $\sF (A)$ of $A$-irreducible subspaces of $\CC^{N}$
is an  orthoframe, and as before    is   invariant under $P$. 
If $B$ is a second such subgroup of $P$ for which 
  $\oA \cap \oB = 0$,   
then $| (u_1,u_2) | = p^{-n/2}=1/\sqrt N$ whenever $u_1$ and $u_2$ are unit vectors in members of $\sF (A) $ and $ \sF (B)$, respectively.

{\em  Each totally isotropic $n$-space} of  $\overline{E}$ arises as
some  $\overline{A}$, and hence {\em  determines a unique orthoframe $\sF (\overline{A}) $.}
  A {\em symplectic spread} of $\overline{E}$ is 
  defined as before.

 Theorem~\ref{Char 2 MUBs} {\em holds with $P$ replaced by $E$\,} \cite[Theorem~11.4, Corollary~11.6]{CCKS}, so that any symplectic spread $\Sigma$ of  $\overline{E}$ produces a
 complete set $ \sF(\Sigma) $ of MUBs  in $\CC^{p^n}$.
Example~\ref{2-dim examples} arises as before, and produces 
  the ``usual'' complete set of MUBs of $\CC^{p^n}$ \cite{GR}; 
  $\Aut  ( \sF(\Sigma))$ behaves   as before.
 Proposition~\ref{MUB test} holds with 2 replaced by $p$.

  The passage from $\Sigma $ to $\sF (\Sigma)$
  is slightly easier to describe in the present setting than 
  in the preceding section.     
First note that (for $p=2$ or $p$ odd) every symplectic spread $\Sigma$ 
  in $\overline{E} = V\oplus V$ can be assumed to be of the following type:
  \begin{equation}
\label{set K} 
 \mbox{\em$\Sigma$ consists of  $\,0\oplus V$ and all $\,\{(v, Mv) \mid v\in V\} $  for  $M\in {\mathcal K},$} \hspace{70pt}
    \end{equation}
 where ${\mathcal K}$ is a set of  $|V|=p^n$ symmetric $n\times n$ matrices such that the difference of any two is nonsingular.
    (This was rediscovered in  \cite[Theorem~4.4]{BBRV},
  without the connection  (\ref{affine plane})  to affine planes.
  The relationship between symplectic spreads and MUBs was also rediscovered in \cite[Sec.~4.5.6]{Howe}, again without the connection to affine planes.)  
  If $p>2$ then  
      \begin{equation}
\label{use spread set}
\begin{array}{llll}
 \mbox{\em$\sF (\Sigma)=\{{\mathcal F}_\infty, \,{\mathcal F}_M^{\mathcal K}  \mid
 M\in {\mathcal K} \},$ where
}\vspace{4pt}
\\
  \mbox{\em${\mathcal F}_\infty :=\! \{\<e_v\>\mid v\in V\}$ and  
${\mathcal F}_M ^{\mathcal K}:= \! \{\<\sum_{v\in V}
\z^{a\cdot v +v\cdot Mv/2}e_v\>\mid a\in V\}$.%
}%
\end{array}%
\hspace{30pt}  
\end{equation}
(Here $v\cdot Mv/2$ is the quadratic form associated with the symmetric bilinear form $u\cdot Mv$.)
The corresponding affine plane has points  $(x,y)\in V\oplus V $ and the following
\begin{equation}
\label{affine plane}
  \mbox{\emph{lines}: \em $x=b\,$   and   $\,y=Mx+b\,$  for   $ \,b\in V,$  $ M \in {\mathcal K}$.} \hspace{110pt}   
      \end{equation}
The simplest ${\mathcal K}$ is $\GF(p^n)$ using 
$M\colon x\mapsto mx$.
Directly verifying  (without use of $P$) that    (\ref{use spread set}) defines MUBs is straightforward:  if $e_{a,M}:=
 \frac{1}{\raisebox{1.7ex} {\hspace{.001pt}}\raisebox{-.2ex} {~}\sqrt N}\sum_{v\in V}
\z^{a\cdot v +v\cdot Mv/2}e_v  $ then 
$(e_{a,M}, e_{a',M'})=
\frac{1}{N}\sum_{v\in V}
\z^{d \cdot v +v\cdot \Delta v/2}   $ with $d:=a-a',$ $ \Delta:=M-M'$,
so that $(e_{a,M}, e_{a,M})=\frac{1}{N}N$.  If $d\ne0$
and $\Delta\ne O$, use  $u=v-v'$ in  the   calculation \vspace{2pt}

$\begin{array}{llllllll}
\vspace{2pt} |(e_{a,M}, e_{a',M'})|^2
\hspace{-6pt}&= & \hspace{-6pt}
\frac{1}{N^2}\sum_{v,v'\in V}
\z^{d \cdot  v  +v\cdot \Delta v/2-d \cdot  v' - v'\cdot \Delta v'/2 } 
\\  
\vspace{2pt}  \hspace{-6pt}&= & \hspace{-6pt}
\frac{1}{N^2}\sum_{u \in V}\z^{d\cdot u -u\cdot \Delta u/2 }
\sum_{v \in V}
\z^{ - v\cdot \Delta u}  
\\  
\vspace{2pt}\hspace{-6pt}&= & \hspace{-6pt}
\frac{1}{N^2} N  ,
\end{array}
$

\vspace{-2pt} \noindent
since $\Delta $   is symmetric and $\Delta u\ne0$ for $ u\ne0$   and $M\ne M'$ in ${\mathcal K}$,
while $\sum_{j=0}^{p-1}\z^j=0 $.

For the case $p=2$ in the preceding section there are minor complications:   the end of (\ref{use spread set})
is replaced by  
\begin{equation}
\label{hats}
{\mathcal F}_M ^{\mathcal K}:= \! \Big\{\<\sum_{v\in V}
i^{2\hat a\cdot\hat  v +\hat v\cdot \hat M\hat v}e_v \> \, \big|\, a\in V \Big\},
\end{equation}
where the ``hats'' denote that the vector or matrix now has entries $0,1$ \emph{viewed inside} $\Z_4$
(so that  $\hat a,\hat  v\in \Z_4^n$).  
Direct verification that we have MUBs is as before with 
additional  bookkeeping.
 The difference between the situations $p=2$ and $p>2$ becomes even more significant when discussing known constructions.

 Most known constructions for odd $p$ are based on generalizations of fields called   \emph{semifields}: algebras satisfying the usual axioms for a field except for the associativity
and commutativity  of  multiplication.
 We refer to \cite{Bier, KaCommutative} for  further information.   Semifields amount to  having  the set  $\mathcal K$ in (\ref{set K})   closed under addition; in this case $\Sigma $ is   called a ``symplectic semifield spread''.
Every commutative semifield corresponds in a somewhat indirect manner to  a symplectic semifield spread (cf. \cite[Proposition~3.8]{KaCommutative}
or the end of Remark~\ref{GR error}; 
 this statement involves two different semifields).
 Therefore we will refer to instances of commutative semifields as if they were examples of symplectic spreads. 
 
  There are   many papers containing (among other things) surveys of commutative semifields, so we mention only two: \cite{Bier, KaCommutative}.  Since there recently have been new constructions for commutative odd order  semifields every few months, 
  the following 
   list is guaranteed to be out of date.
   
\begin{Example} \rm 
\label{odd examples}
We survey   the known examples of symplectic spreads of  $\overline{E}\cong \Z_p^{2n}$ for odd $p$, and hence implicitly the corresponding complete sets of MUBs in~$\CC^{p^n}\!$
obtained as in  Theorem~\ref{Char 2 MUBs}. 
We give   explicit complete sets of MUBs in (b).
\begin{itemize}
\item[(a)]  The analogue of Example~\ref{2-dim examples} 
uses a 2-dimensional vector space over $\GF(p^n)$.
\label{odd 2-dim examples}
 {\em These are the only complete sets of} MUBs  {\em arising in} Proposition~\ref{MUB test}  {\em when $N=p^n=p$}, i.\hspace{1pt}e., when $|E|=p^3$.
\item[(b)] Older families \cite{Dic, Alb} have an unbounded    number of pairwise inequivalent examples     as a function of $n$.
  \emph{We present examples
  not unitarily equivalent to the  one in} (a) \emph{nor to one another}, based on \cite{BKL}. 
  Let $V$ be   $K=\GF(p^n)$ with $n$ odd, choose an integer $s$  relatively prime to $n$ such that $1\le s <n/2$, 
  and let  $T(x):=\sum_{j=0}^{n-1}x^{p^j}$, ${x\in K}$, be the trace map.  
Label the standard orthonormal basis  of $\CC^{p^n}$ as $e_x,\,  x\in K$,
with corresponding orthoframe  $\sF_\infty$. Then $\sF :=  \{\sF_\infty, \, \sF[b]\mid b\in K\}$ 
\emph{is a complete set of} MUBs, where 
$$ \qquad \quad
\sF[b]:=\Big\{\<\sum_{x\in K}\z^{T(a x) \, +\, 
T(bx^{p^{n-s}+1}  + \, b^{p^s} x^{p^s+1} ) / 2}e_x\> \, \big|\,  a\in K \Big \}.
\vspace{-4pt}
$$ 
If  the exponents are written $ {a \cdot x }+
 {x \cdot ({b x^{p^{n-s}}}\! + {b^{p^s} x^{p^s}})/2}$ 
using a dot  product on $K$ as in (\ref{use spread set}),
 the result is an equivalent set of MUBs.
Allowing $s=0$ would give  (a). 
\item[(c)] Families  with $p=3$ 
\cite{CG,Ga, TP, PW, CM,DiY}:    for some of these  the 
  number of pairwise inequivalent examples is unbounded   as a function of $n$.
\item[(d)] Recent families 
\cite{ZW, BuHe, LMPT, Bier, BK}:  for some of these  the 
  number of pairwise inequivalent examples is unbounded   as a function of $n$.
\item[(e)] There are also two families of examples not related to semifields
\cite{KaOvoids, BBLP}, the first having $p=3$. 
  This is very different from the situation in characteristic 2, where Examples~\ref{symplectic char 2}(b)-(d) are fairly  large families not   corresponding to semifields.
(Examples~\ref{symplectic char 2}(a)   correspond to semifields.)
\end{itemize}

The number of items in the above list attests to the amount of research occurring on this topic.  Nevertheless, 
each of these families is associated with  fewer than   $\sqrt N$ pairwise inequivalent  complete sets of MUBs in $\CC^N$.
In view of Example~\ref{symplectic char 2}(a), this means that at present there are far fewer complete sets that have been obtained using odd characteristic than there are using characteristic 2.
As in the preceding section,  there are undoubtedly large numbers of examples yet to be found.
\end{Example}

 \begin{Remark}
 \label{GR error}
 \rm
  Contrary to \cite[p.~255]{GR}, \emph{it is not the case that} ``there is a natural correspondence between semifields and symplectic spreads (see
for example'' \cite [Proposition~3.8]{KaCommutative}). The cited result in \cite{KaCommutative} 
 only concerns commutative semifields and symplectic 
 {\em semifield} spreads, not general  symplectic spreads.
Therefore, contrary to \cite[p.~255]{GR}, 
it is not the case that
their ``mutually unbiased bases \dots~ are
equivalent to those of'' \cite{CCKS}. 
Examples~\ref{symplectic char 2}(b)-(e)
and \ref{odd examples}(e)
do not arise from  semifields other than fields.

On the other hand,  the complete sets of MUBs obtained   
in \cite{GR} (by using 
 $ \sum_{v\in V}
\z^{a\cdot v +b\cdot (v*v)/2}$ in place of
the vector sum in    (\ref{use spread set})
for a
commutative semifield multiplication $*$  on $V$) are seen to be among those in \cite{CCKS} by a straightforward use of Proposition~\ref{MUB test}.
Even simpler is to match up with (\ref{use spread set}):
if  we let $\e_i$  denote the $i$th  standard basis vector of $V$,  then
 ``solve'' $b\cdot (x*y)=x\cdot M y$ for the symmetric matrix $M=(M_{ij}) $ as a function of $b\in V$ via
$M_{ij} =\e_i\cdot M \e_j=b \cdot (\e_i*\e_j)$,
so that $b\cdot (v*v)=v\cdot Mv$. 
 \end{Remark}

  The complete
sets of MUBs described in \cite{GR} are also obtained in \cite{RS}, but
the latter paper goes further:

\begin{Example} \rm 
\label{planar function examples}
Complete sets   $ \sF ^f$  of MUBs in $\CC^{3^n}\!$ 
are obtained in \cite{RS}
for each odd $n\ge 5$, corresponding to the ``planar  
functions''\footnote{This means that $f(x+a)-f(x)=b$ has a unique solution $x$  for any  $a\ne0$ and $b$ in $K$.}
 $f(x)=x^{(3^k+1)/2}$ on $K=\GF(3^n)\cong V$ (where $2n$ and $k$ are relatively prime and $k\not \equiv \pm1$ (mod $2n$))   
  \cite[Theorem~6.2]{CM}. This time  ${\mathcal F} ^f$  consists of 
${\mathcal F}_\infty$  and all 
${\mathcal F}_b^f:=
\{\< \sum_{v\in K}\z^{a\cdot v+b\cdot f(v)}e_v \>\mid a\in K \}$ for $b\in K$.  A direct proof that $\sF ^f$ is a set of MUBs is similar to the one following (\ref{affine plane}),   using 
$v-v'=u$ and 
$\z^{(a-a')\cdot v+ (b-b')\cdot f(v)-(a-a')\cdot v' - (b-b')\cdot f(v')}  =
\z^{(a-a')\cdot u}\z^{ (b-b')\cdot (f(v)-f(v-u))}$.
 It is easy to
use Proposition~\ref{MUB test}  to show that these  do not arise from 
symplectic spreads in   any extraspecial $p$-group
(this equivalence  question is not discussed in \cite{RS}).  {\em These are the only known complete sets  of} MUBs  {\em that do not arise from 
symplectic spreads using   extraspecial groups.}
These have a property in common with complete sets arising 
from symplectic semifields:  there is a group of $N^2$ automorphisms 
(generated by all $e_x \mapsto \z^{c\cdot x}  e_x$ and all $e_x\mapsto \z^{c\cdot f(x)}  e_x$, $c\in V$)
having orbits of size $1$ and $N^2$ on
 $\cup\sF ^f$.

This unusual family of complete sets of MUBs suggests that there are many more families yet to be discovered that are not related to extraspecial groups.
 
\end{Example}

 \begin{Remark}
 \label{Where are the planes}\bf Where are the planes?
 \rm Each of the known examples in Sections~\ref{MUBS characteristic 2}  and \ref{MUBS  $p>2$} has at least one associated affine plane. 
 \emph{Where are these  planes}?
 
 For complete sets $\sF$ of MUBs obtained as in Theorem~\ref{Char 2 MUBs}, the answer is similar to Proposition~\ref{MUB test}:  $\sF$ uniquely determines a group $P$ or $E$ of automorphisms, hence also a symplectic spread 
 $\Sigma$ and set ${\mathcal K}$
 as in 
 (\ref{set K}), and finally an affine plane $\pi$ as in 
 (\ref{affine plane}).
Moreover, by  (\ref{set K}) and (\ref{affine plane}), $\sF$ can be identified with the set $\Sigma$ 
of parallel classes of  $\pi$,  
and then points can be identified with
$V\oplus V,$ or equivalently, with (suitable!)   subsets of  $\cup\sF$  consisting of one  element  from each member  of $\sF$.
 
 For complete sets obtained as in     Example~\ref{planar function examples}, 
the associated affine  plane $\pi(f)$  has as points the vectors in $V\oplus V$ and as lines
the sets  $x=b$ 
and $y=f(x+a)+b$ \ ($a,b\in K$).  It would be helpful to ``see'' this plane
in terms of $\sF ^f$, perhaps  
using  $\sF ^f$ as the set of $N +1$ parallel classes of lines and 
(as above) ``natural'' subsets of  $\cup\sF ^f$ 
corresponding to  points.
However, \emph{such a  description could not be invariant under the 
automorphism 
group of order $N^2$ mentioned in} Example~\ref{planar function examples}:  $\pi(f)$   has no automorphism group of that order inducing the identity on one of its parallel classes.\footnote{Curiously, there is also a group of $N^2$ automorphisms of $\pi(f)$  that does not act on $\sF ^f$.  (This is the group of all $(x,y)\mapsto (x+c,y+d)$, $c,d\in K$, having orbits of size $N$ and $N^2$ on the set of all lines.)} 
 In other words, there is no canonical 
 (i.\:e., $\Aut(\sF ^f$)-invariant) way to 
 obtain such a description from $\sF ^f$.   Nevertheless, it is at least somewhat plausible that there might be 
 an entirely different way to reconstruct 
 $\pi(f)$ either  from~$\sF^f$ or from some entirely different complete set of MUBs associated with $f$.
 \end{Remark}

\vspace{-8pt}

\section{Complete sets of MUBS in $\R^{2^n}$} 
\label{Real MUBs}

 We return to the group $E$ in Section~\ref{MUBS characteristic 2}. 
 This time we restrict to vectors in $\R^N,N=2^n$, and use the usual inner product  $(~,~)$.
  There is additional structure to consider:  the function $Q\colon E \to \Z_2$ given by  $Q(x)=x^2\in \zen(E)$ determines a quadratic form
   $\overline Q\colon  \overline{E}\to \Z_2$
   that polarizes to 
 the alternating bilinear form {\bf (~,~)} appearing in  Section~\ref{MUBS characteristic 2}   (i.e., 
 $\overline Q(x+y)-\overline Q(x)-\overline Q(y)$=\mbox{\bf ($x , y$)} for all $x,y\in \overline{E}$).
 This time we are interested in subgroups $A$ of $E$ such that 
$\overline{A}$ is {\em totally singular}:  $\overline Q(\overline{A})=0$ (and hence also {\bf ($\overline{A} , \overline{A}$)} $=0)$.  If $\overline{A}$ is a  totally singular $n$-space then the set
$\sF (A)=\sF (\overline A)$ of $A$-irreducible subspaces of $\R^{N}$
is an {\em orthoframe}: a set of $N$ pairwise
orthogonal 1-spaces. Once again, if $B$ is a second such subgroup of 
$E$ for which 
  $\oA \cap \oB = 0$,   
then $| (u_1,u_2) | =  1/\sqrt N$ whenever $u_1$ and $u_2$ are unit vectors in members of $\sF (A) $ and $ \sF (B)$, respectively. 
Any family $\sF$ of orthoframes satisfying this last property involving pairs of unit vectors has size at most $\frac{1}{2} N  +1 $ \   \cite[(3.9)]{CCKS}. \vspace{2pt}
When equality holds we have a \emph{complete set of} MUBs {\em of $\R^N$.} Note that such  sets are smaller than the complete sets in the previous sections. Moreover,  the factor  $\frac{1}{2}$ leads us only  to  use vector spaces $\overline{E}$ of characteristic 2.
See \cite{BSTW1} for more information concerning sets of real MUBs. 
  
 An {\em orthogonal spread} of $\overline{E}$ is a family $\Sigma$ of
 $ 2^{n-1}+1$ totally singular $n$-spaces 
 of $\overline{E}$ 
 any two of which have zero intersection; every member of 
 $\{ 0\ne x\in \overline{E}\mid \overline Q(x)=0 \}$ 
 lies in one and only one member of $\Sigma$.  
 As in Theorem~\ref{Char 2 MUBs} there is a corresponding complete set  $\sF (\Sigma) $  of MUBs of $\R^N$ \cite[Theorem~3.4]{CCKS}; and the equivalence part of the theorem continues to hold   
 \cite[Proposition~3.16]{CCKS}.
However, each orthogonal spread is associated with a somewhat large number of possibly nonisomorphic affine planes \cite[Sec.~7]{CCKS}.

\begin{Example} \rm 
All orthogonal  spreads of   $\Z_2^{2n}$  arise from symplectic spreads of $\Z_2^{2n-2}$, and $n$ {\em must be even}
 (\cite{Di,Dy}  or \cite[I,~Sec.~3]{KaKer1}). One of these orthogonal  spreads corresponds to 
Example~\ref{2-dim examples}; the other  known ones correspond to     Examples~\ref{symplectic char 2}(a)-(d)
(note that $n$ here corresponds to $n-1$ in those examples). All  comments about numbers of inequivalent complete sets of MUBs   hold
as before.

In fact the relationship with orthogonal spreads is an essential ingredient for all of  the constructions in Examples~\ref{symplectic char 2}(a)-(d). 


\end{Example}

{\noindent \bf Acknowlegement.}  I am grateful to Mary Beth Ruskai for many very helpful comments.

\end{document}